\documentclass{article}

\usepackage[utf8]{inputenc}
\usepackage[english]{babel}

  \usepackage[T1]{fontenc}
  \usepackage[mathcal]{eucal}
  \usepackage{amsfonts}
 \usepackage[colorlinks=false]{hyperref}
  
  \usepackage{amsmath}
  \usepackage{amsthm}
  \usepackage{amssymb}
  \usepackage{mathrsfs}
  \usepackage{stmaryrd}
  \usepackage{graphicx}
  \usepackage[all]{xy}
  \usepackage{relsize, exscale} 
  \usepackage{tikz-cd}
  \usepackage{wasysym}
  \usepackage{enumitem}
  \usepackage{float}
\usepackage{authblk}
\usepackage{hyperref}
\usepackage{mathtools}
  
  \newsavebox\mybox
\newlength\mylength
\newlength\mylengthb

\newtheorem{thm}{Theorem}[section]
\newtheorem{df}[thm]{Definition}
\newtheorem{thdf}[thm]{Theorem-Definition}
\newtheorem{lem}[thm]{Lemma}
\newtheorem{prop}[thm]{Proposition}
\newtheorem{cor}[thm]{Corollary}
\newtheorem{nota}[thm]{Notation}
\newtheorem{rem}[thm]{Remark}
\newtheorem{ex}[thm]{Example}

\newtheorem{set}[thm]{Setup}

\newcommand{\bea}{\begin{eqnarray}}
\newcommand{\eea}{\end{eqnarray}}
\newcommand{\bna}{\begin{eqnarray*}}
\newcommand{\ena}{\end{eqnarray*}}

\newcommand{\Op}{\mathcal{O}_{\mathbb{P}^1}}
\newcommand{\Oe}{\mathcal{O}_{\mathbb{P}(\mathcal{E})}}

\newcommand{\Pic}{\operatorname{Pic}}
\newcommand{\Bl}{\operatorname{Bl}}
\newcommand{\NS}{\operatorname{NS}}

\newcommand{\Bs}{\operatorname{Bs}}

\newcommand{\Ext}{\operatorname{Ext}}
\newcommand{\Sym}{\operatorname{Sym}}
\newcommand{\dg}{\operatorname{deg}}

\makeatletter
\newcommand{\subjclass}[2][1991]{%
  \let\@oldtitle\@title%
  \gdef\@title{\@oldtitle\footnotetext{#1 \emph{Mathematics subject classification.} #2}}%
}
\newcommand{\keywords}[1]{%
  \let\@@oldtitle\@title%
  \gdef\@title{\@@oldtitle\footnotetext{\emph{Key words and phrases.} #1.}}%
}

\title{\textbf{Rationally connected threefolds with nef and bad anticanonical divisor}}
\author{Zhixin Xie}
\date{}
\begin{document}
\maketitle
\begin{abstract}
Let $X$ be a smooth projective rationally connected threefold with nef anticanonical divisor. We give a classification for the case when $-K_X$ is not semi-ample.
\end{abstract}

\section{Introduction}
Let $X$ be a complex projective manifold. We say that $X$ is a Fano manifold if the anticanonical divisor $-K_X$ is ample. The classification of three-dimensional Fano manifolds by Mori and Mukai (for $\rho>1$, \cite{MR641971}) and by Iskovskih (for $\rho=1$, \cite{MR463151} \cite{MR503430}) is one of the first achievements of the minimal model program with an impressive number of applications. Projective manifolds with nef anticanonical divisor $-K_X$ are a natural generalisation of Fano manifolds, and one hopes to similarly fulfil a complete classification for this class of manifolds. In \cite{MR3959071}, Cao and H\"{o}ring
showed a decomposition theorem for these manifolds: the universal cover $\tilde{X}$ of $X$ decomposes as a product
\[
\tilde{X}\simeq \mathbb{C}^q\times\prod Y_j\times\prod S_k\times Z,
\]
where $Y_j$ are irreducible Calabi-Yau manifolds, $S_k$ are irreducible hyperk\"{a}hler manifolds, and $Z$ is a rationally connected manifold.

In view of this result, the most interesting case is when $X$ is rationally connected: it is also the most difficult one. Before considering classification for this class of varieties, one should check whether in any given dimension there are only finitely many families of varieties for this class. This property is known under the name of boundedness (see Definition \ref{dfbounded}). Recently, Birkar, Di Cerbo and Svaldi proved in \cite[Theorem 1.6]{birkar2020boundedness} that birationally, there are only finitely many deformation families of projective rationally connected threefolds with $\epsilon$-lc singularities and nef but not numerically trivial anticanonical divisor. Thus it is in principle possible to classify these varieties as has been done for Fano threefolds. If the anticanonical divisor is semi-ample, then there is a standard approach to the classification via a discussion of the anticanonical map and extremal contractions.

We call a nef divisor good if its Iitaka dimension and its numerical dimension are equal, otherwise we call it bad. By a result of Kawamata \cite[Theorem 6.1]{MR782236}, if $-K_X$ is good then it is semi-ample. In this paper we focus on the much more delicate case where $-K_X$ is nef but not semi-ample. Bauer and Peternell have shown in \cite[Theorem 1]{MR2129540} that this implies that the nef dimension (see Definition \ref{nefdim}) $n(-K_X)=3$, the numerical dimension $\nu(-K_X) = 2$ and the Iitaka dimension $\kappa(-K_X)=1$ (in fact they show that $h^0(X,\mathcal{O}_X(-K_X))\geq 3$). It is thus natural to start the investigation with the base locus of the anticanonical system. We start by addressing the case in which the base locus has no divisorial part - a case that was not covered in \cite[Proposition 7.2]{MR2129540}. The first main theorem of this paper is the following effective classification result for this case and we will give some examples in the first part of the paper (Section \ref{3}):
\begin{thm}\label{mainmobile}
Let $X$ be a smooth projective rationally connected threefold $X$ with $-K_X$ nef, $n(-K_X)=3$, $\nu(-K_X)=2$. Suppose that the anticanonical system has no fixed part. Then $-K_X\sim 2D$ where $D$ is a Cartier divisor, and $X$ is one of the following:
\begin{enumerate}[label=\normalfont(\arabic*)]
    \item $X\to \mathbb{P}^1$ is a del Pezzo fibration with general fiber isomorphic to $\mathbb{P}^1\times\mathbb{P}^1$. Then $X\subset\mathbb{P}_{\mathbb{P}^1}(\mathcal{E})$ with
    \[
\mathcal{E}=\mathcal{O}_{\mathbb{P}^1}\oplus\mathcal{O}_{\mathbb{P}^1}\oplus\mathcal{O}_{\mathbb{P}^1}(-1)\oplus\mathcal{O}_{\mathbb{P}^1}(-1),
\]
    and $X$ is an element of the linear system $|\mathcal{O}_{\mathbb{P}(\mathcal{E})}(2)+4F|$, where $F$ is a general fiber of $\pi\colon \mathbb{P}_{\mathbb{P}^1}(\mathcal{E})\to \mathbb{P}^1$.
\item $X=\mathbb{P}(\mathcal{E})$ is a $\mathbb{P}^1$-bundle over a smooth rational surface $Y$ with $-K_Y$ nef, where $\mathcal{E}$ is a nef rank $2$ vector bundle with $c_1(\mathcal{E})=-K_Y$ and $c_2(\mathcal{E})=K_Y^2$, given by an extension
$$
0\to\mathcal{O}_Y\to \mathcal{E}\to\mathcal{I}_Z\otimes \mathcal{O}_Y(-K_Y)\to 0,
$$
where $\mathcal{I}_Z$ is the ideal sheaf of $c_2(\mathcal{E})$ points on $Y$.
\item $X=\Bl_p(Y)$ is the blow-up in a point $p$ of a smooth almost del Pezzo threefold $Y$ of degree $1$ such that $p$ is not the base point of $|-\frac{1}{2}K_Y|$.
\end{enumerate}

Conversely, let $X$ be a variety that appears in one of the above cases with respectively the following conditions:
\begin{enumerate}[label=\normalfont(\arabic*')]
    \item $X\in|\mathcal{O}_{\mathbb{P}(\mathcal{E})}(2)+4F|$ is a very general member;
\item the finite set $Z$ of $c_2(\mathcal{E})$ points are in sufficiently general position on $Y$;
\item the blown up point $p\in Y$ is sufficiently general.
\end{enumerate}
Then $-K_X$ nef, not semi-ample and $-K_X\sim 2D$ where $D$ is a Cartier divisor.
\end{thm}
The class of varieties in case $(3)$ of the above theorem also appeared in \cite[Section 2]{MR3510909} where an explicit and pathological example is constructed. For a complete classification list of smooth almost del Pezzo threefolds of degree $1$, we refer to \cite{MR2427467}. Hence the smooth threefolds with nef and not semi-ample anticanonical divisor whose base locus has no divisorial part are completely classified. 

The second part of the paper (Section \ref{4}) deals with the case when the base locus has a divisorial part. It turns out, that after a sequence of flops, one can assume that the mobile part is always nef. We will show that the mobile part does not have base locus of codimension $2$ and to be more precise:
\begin{thm}\label{mainfixed}
Let $X$ be a smooth projective rationally connected threefold $X$ with $-K_X$ nef, $n(-K_X)=3$, $\nu(-K_X)=2$. Suppose that the anticanonical system has non-empty fixed part. Then there exists a finite sequence of flops $\psi\colon X\dashrightarrow X'$ such that the following holds: 
\begin{itemize}
\item $X'$ is smooth,
    \item $-K_{X'}$ is nef,
    \item the mobile part $|B'|$ of $|-K_{X'}|$ is nef.
\end{itemize}
In this case, $B'^2=0$ and $|B'|$ is base-point-free. It induces a fibration $f\colon X'\to\mathbb{P}^1$.
\end{thm}

Back to the problem of birational boundedness for the family of smooth projective rationally connected threefolds with nef (and not semi-ample) anticanonical divisor, we deduce the boundedness from Theorem \ref{mainmobile} (and Proposition \ref{conicbundle}) for the case when the anticanonical system has no fixed part. As for the case when there is a non-empty fixed part, the boundedness does not follow from Theorem \ref{mainfixed}. However we expect a complete classification for this case as in the previous one. In the case of Theorem \ref{mainfixed}, the structure of $X'$ can be complicated. Examples of such $X'$ which are not isomorphic to a product can be found when the general fiber of $f$ is $\mathbb{P}^2$ blown up in $9$ points such that the unique element in the anticanonical system is a smooth elliptic curve, we refer to the author's thesis \cite{Xie2021}.

We give now a sketch of proof of Theorem \ref{mainfixed}. The idea is to argue by contradiction and suppose that the mobile part of the anticanonical system has base locus of codimension two. We first give a description of the anticanonical system in this case:
\begin{prop}\label{descrip}
Let $X$ be a smooth projective rationally connected threefold $X$ with $-K_X$ nef, $n(-K_X)=3$, $\nu(-K_X)=2$. Suppose that the anticanonical system has non-empty fixed part and that the mobile part $|B|$ is nef. If $B^2\neq 0$, then $|-K_{X}|=A+|2H|$ where $H$ is a prime divisor. In this case, let $F$ be a general member in $|H|$. Then both $A$ and $F$ are $\mathbb{P}^1$-bundles over a smooth elliptic curve, such that their anticanonical divisors $-K_A$ (resp. $-K_F$) are nef and divisible by $2$ in $\Pic(A)$ (resp. $\NS(F)$). Furthermore, both $A\cdot F$ and $F^2$ are smooth elliptic curves.
\end{prop}
By running the minimal model program, with the classification of contractions by Mori for smooth threefolds \cite{MR662120}, we show in this case:
\begin{prop}\label{mmp}
In the setting of Proposition \ref{descrip}, 
there exists a finite sequence \[ X=X_0\overset{\varphi_1}\rightarrow X_1\overset{\varphi_2}\rightarrow\cdots \overset{\varphi_k}\rightarrow X_k\overset{\varphi_{k+1}}\rightarrow Y \] with $k=0$ or $1$, where $X_i$ is a smooth threefold with $-K_{X_i}$ nef such that $|-K_{X_i}|$ has non-empty fixed part, $\varphi_i$ is a blow-up along a smooth elliptic curve and
$Y$ is a smooth threefold with $-K_Y$ nef, $n(-K_Y)=3$, $\nu(-K_Y)=2$ such that $|-K_Y|=|2G|$ has no fixed part and a general member in $|G|$ is isomorphic to $F$.
\end{prop}
The threefold $Y$ is classified in Theorem \ref{mainmobile}. Comparing the general member of $|-\frac{1}{2}K_Y|$ in Theorem \ref{mainmobile} with the geometry of the surface $F$ will then lead to a contradiction.

{\bf Acknowledgements.} I would like to express my deep gratitude to my supervisor Professor Andreas H\"{o}ring for numerous discussions and constant support. I also heartily thank Professor Cinzia Casagrande for her careful proofreading and valuable comments.  I thank the referee for suggestions to improve the exposition. The author is supported by the IDEX UCA JEDI project (ANR-15-IDEX-01) and the MathIT project.
\section{Preliminaries}

\begin{df}\label{dfbounded}
A set of projective varieties $\mathbb{X}$ is said to be bounded if there exists $\phi\colon\mathfrak{X}\to\mathfrak{B}$ a projective morphism of schemes of finite type such that each $X\in\mathbb{X}$ is isomorphic to $\mathfrak{X}_b$ for some closed point $b\in\mathfrak{B}$.
\end{df}

\begin{nota}\textit{\cite[Definitiion 2.1.3, Remark 2.3.17]{MR2095471}
Let $X$ be a normal projective variety and $D$ an $\mathbb{R}$-Cartier divisor on $X$. We denote by
\begin{itemize}
    \item $\kappa(D)$ the Iitaka (Kodaira) dimension of $D$.
    \item $\nu(D)\coloneqq {\rm{max}}\ \{n\mid D^n\not\equiv 0\}$ the numerical dimension of $D$ when $D$ is nef.
\end{itemize}
}
\end{nota}

\begin{thdf}\cite[Theorem 2.1]{MR1922095}\label{nefdim}
Let $L$ be a nef line bundle on a normal projective variety $X$. Then there exists an almost holomorphic dominant meromorphic map $f\colon X\dashrightarrow B$ with connected fibers such that
\begin{enumerate}[label=\normalfont(\arabic*)]
    \item $L$ is numerically trivial on all compact fibers $F$ of $f$ of dimension ${\rm{dim}}\ X - {\rm{dim}}\ B$
    \item for a general point $x\in X$ and every irreducible curve $C$ passing through $x$ such that ${\rm{dim}}\ f(C) > 0$, we have $L\cdot C > 0$.
\end{enumerate}
The map $f$ is unique up to birational equivalence of $B$. In particular ${\rm{dim}}\ B$ is an invariant of $L$ and we set $n(L) = {\rm{dim}}\ B$, the nef dimension of $L$.
\end{thdf}

Note that the nef dimension of a nef line bundle $L$ on $X$ is maximal, i.e. $n(L)=\dim X$, if and only if the variety $X$ is not covered by $L$-trivial curves. Moreover, we have the following inequalities which relate the above three invariants:
\begin{thm}\cite[Proposition 2.2]{MR782236}, \cite[Proposition 2.8]{MR1922095}\label{inequality}
Let $X$ be a smooth projective variety. Let $L$ be a nef divisor on $X$. Then we have
\[
\kappa(L)\leq \nu(L)\leq n(L).
\]
\end{thm}

\begin{lem}\label{nef}
Let $X$ be a normal projective $\mathbb{Q}$-factorial variety with $-K_X$ nef. Let $D$ be an effective $\mathbb{Q}$-divisor such that the pair $(X,D)$ is log canonical. If $D$ is not nef, then there exists a $(K_X+D)$-negative extremal ray $\Gamma$ such that $D\cdot\Gamma<0$.
\end{lem}
\begin{proof}
 Suppose that there is no such extremal ray. Since $D$ is not nef, there exists an irreducible curve $l\subset X$ such that $D\cdot l<0$. Then we can write
\[
l=\sum_i \lambda_i\Gamma_i + R,
\]
where 
\begin{itemize}
    \item $\lambda_i\geq 0$;
    \item $\Gamma_i$ are $(K_X+D)$-negative extremal rays. By assumption they all satisfy $D\cdot \Gamma_i\geq 0$;
    \item $(K_X+D)\cdot R\geq 0$.
\end{itemize}
Therefore,
\[
0>D\cdot l = \sum_i \lambda_i D\cdot\Gamma_i + D\cdot R\geq D\cdot R,
\]
i.e. $D\cdot R<0$.

Since $(K_X+D)\cdot R\geq 0$, we have
\[
K_X\cdot R\geq -D\cdot R >0,
\]
which contradicts the fact that $-K_X$ is nef.
\end{proof}

\begin{lem}\label{mobilepart}
Let $X$ be a smooth projective rationally connected threefold. Let $D$ be a divisor with $\kappa(D)=1$. Suppose that the linear system $|D|$ has no fixed part and the general member in $|D|$ is reducible. Then a general member in $|D|$ is linearly equivalent to $mH$, where $H$ is a prime divisor and $m\geq 2$. Furthermore, $h^0(X,\mathcal{O}_X(H))=2$ and $h^0(H,\mathcal{O}_H(H|_H))=1$.
\end{lem}
\begin{proof}
Let $\phi\colon X\dashrightarrow C$ be the rational map determined by the linear system $|D|$. Then $C\simeq \mathbb{P}^1$ as $\kappa(D)=1$ and $h^1(X,\mathcal{O}_X)=0$.

Let $\mu\colon\tilde{X}\to X$ be a birational modification which resolves the base locus of $|D|$. Let $F$ be a general fiber of the induced morphism $\tilde{\phi}\colon\tilde{X}\to C$. Since $|D|$ has no fixed part, the pushforward $\mu_*(F)$ is a general member of $|D|$. Furthermore, the general fiber $F$ is not connected as the general member in $|D|$ is reducible.

Let $\tilde{\phi}'\colon\tilde{X}\to C'$ be the Stein factorization of the morphism $\tilde{\phi}$ and $\nu\colon C'\to C$. Then $\tilde{\phi}'$ has smooth connected general fiber and $C'\simeq \mathbb{P}^1$ as $h^1(X,\mathcal{O}_X)=0$.

For a point $p\in C$, we have 
\[
\nu^*(p)\simeq \mathcal{O}_{\mathbb{P}^1}(m)
\]
with $m\geq 2$ the number of connected components of $\tilde{\phi}^*(p)$ and thus \[
\tilde{\phi}^*(p) = \tilde{\phi}'^*(\nu^*(p))\simeq \tilde{\phi}'^*(\mathcal{O}_{\mathbb{P}^1}(m)).
\]
Let $F'$ be a general fiber of $\tilde{\phi}'\colon\tilde{X}\to\mathbb{P}^1$. Then $F\sim mF'$ and thus a general member in $|D|$ is linearly equivalent to $mH$ where $H\coloneqq\mu_*(F')$. Hence \[h^0(X,\mathcal{O}_X(H))=h^0(\mathbb{P}^1,\mathcal{O}_{\mathbb{P}^1}(1))=2.\]
Now the exact sequence:
\[
0\to \mathcal{O}_X\to\mathcal{O}_X(H)\to \mathcal{O}_H(H|_H)\to 0
\]
gives $h^0(H,\mathcal{O}_H(H|_H))=1.$
\end{proof}

From now onwards, $X$ will denote a smooth projective rationally connected threefold with $-K_X$ nef, $n(-K_X)=3$ and $\nu(-K_X)=2$. 

By \cite[Theorem 2.1]{MR2129540}, the condition of $-K_X$ satisfying $n(-K_X)=3$ and $\nu(-K_X)=2$ is equivalent to $\nu(-K_X)=2$ and $\kappa(-K_X)=1$. The latter one is more useful since it is in practice easier to compute the Iitaka dimension than the nef dimension.

Since $X$ is rationally connected, we have $\chi(\mathcal{O}_X)=1$. Together with $K_X^3 = 0$, we deduce by Riemann-Roch theorem that $\chi(-K_X) = 3$. Moreover, as $-K_X$ is nef and $\nu(-K_X)=2$, by Kawamata-Viehweg vanishing theorem \cite[Corollary]{MR675204}, one has $H^1(X, \mathcal{O}_X(2K_X))=0$. Hence we deduce that $$H^2(X,\mathcal{O}_X(-K_X)) = 0.$$ by Serre duality. Therefore,
\[
h^0(X,\mathcal{O}_X(-K_X))\geq 3.
\]

\begin{lem}\label{2irred}
Let $X$ be a smooth projective rationally connected threefold with $-K_X$ nef, $n(-K_X)=3$ and $\nu(-K_X)=2$. Let $|B|$ be the mobile part of the anticanonical system $|-K_X|$ and $D$ be a general member in $|B|$. Then $D$ has at least two irreducible components. 
\end{lem}
\begin{proof}
We can write $|-K_X|=A+|B|$ with $A$ the fixed part (which can be empty) and $|B|$ the mobile part. For a general member $D$ of $|B|$, we have the following exact sequence:
\[
0\to \mathcal{O}_X(-K_X-D)\to \mathcal{O}_X(-K_X)\to \mathcal{O}_D(-K_X|_D)\to 0.
\]
Since \[h^0(X,\mathcal{O}_X(-K_X-D))=h^0(X,\mathcal{O}_X(A))=1,\] together with $h^0(X,\mathcal{O}_X(-K_X))\geq 3,$ we have
$h^0(D,\mathcal{O}_D(-K_X|_D))\geq 2.$

Now suppose by contradiction that $D$ is irreducible.
Let $\nu\colon\overline{D}\to D$ be the normalization of the surface $D$. Then for the pullback of the Cartier divisor $-K_X|_D$, we have
\[
h^0(\overline{D},\nu^*(-K_X|_D))\geq h^0(D,-K_X|_D)\geq 2.
\]
Hence the linear system $|\nu^*(-K_X|_D)|$ on $\overline{D}$ has a mobile part $M$. On the other hand,
since $-K_X$ is nef and $(-K_X)^3=0$, one has $(-K_X)^2\cdot D=0$, i.e. $(-K_X|_D)^2=0$.
Since $\nu^*(-K_X|_D)$ is nef and $\nu^*(-K_X|_D)^2=(-K_X|_D)^2=0$, we deduce that \[\nu^*(-K_X|_D)\cdot M=0.\] Therefore, $\overline{D}$ is covered by $\nu^*(-K_X|_D)$-trivial curves, from which we deduce that $D$ is covered by $(-K_X)$-trivial curves. As $D$ moves, this contradicts the fact that $n(-K_X)=3$.
\end{proof}

Now Lemmas \ref{mobilepart} and \ref{2irred} give the following:
\begin{cor}\label{mH}
Let $X$ be a smooth projective rationally connected threefold with $-K_X$ nef, $n(-K_X)=3$ and $\nu(-K_X)=2$. Let $|B|$ be the mobile part of the anticanonical system $|-K_X|$. Then
\[
B\sim mH,
\]
where $m\geq 2$ and $H$ is a prime divisor such that $h^0(H,\mathcal{O}_H(-K_X|_H))=1$, $h^0(X,\mathcal{O}_X(H))=2$ and $h^0(H,\mathcal{O}_H(H|_H))=1$.
\end{cor}
\begin{proof}
It remains to show that
$
h^0(H,\mathcal{O}_H(-K_X|_H))=1.
$
By contradiction, suppose that $h^0(H,\mathcal{O}_H(-K_X|_H))\geq 2$. In Lemma \ref{2irred}, we may repeat the same argument in the second part of the proof with $H$ playing the role of $D$, then the argument following from the normalization of the surface leads to a contradiction.
\end{proof}

\begin{lem}\label{lemconverse}
Let $X$ be a smooth projective threefold with $-K_X$ non-zero effective, divisible by two in $\Pic(X)$ and $K_X^3=0$. Suppose that there exists an irreducible normal surface $H\in|-\frac{1}{2}K_X|$ such that $-K_H$ is nef, non-zero effective and not semi-ample. Then $-K_X$ is nef and not semi-ample, i.e. $\nu(-K_X)=2$ and $\kappa(-K_X)=1$.
\end{lem}
\begin{proof}
We have $-K_X\sim 2H$. The adjunction formula gives $-K_H\sim H|_H$. 
We first show that $H$ and thus $-K_X$ is nef.
Indeed, it is enough to show that the restriction of $H$ on itself is nef: let $C\subset H$ be an integral curve, then \[ H\cdot C=H|_H\cdot C=-K_H\cdot C\geq 0 \]
as $-K_H$ is nef.

Now since there exists a non-zero effective divisor in $|-K_H|$, we deduce that $\nu(-K_X)=\nu(H)=2$.

Since $-K_H$ is not semi-ample, we have $\kappa(H,-K_H)=0$. Then for any $m\geq 1$, \[ h^0(H,\mathcal{O}_H(mH|_H))=h^0(H,\mathcal{O}_H(-mK_H))=1.\]
Now the short exact sequence
\[
0\to\mathcal{O}_X((m-1)H)\to\mathcal{O}_X(mH)\to \mathcal{O}_H(mH|_H)\to 0
\]
gives $h^0(X,\mathcal{O}_X(mH))\leq h^0(X,\mathcal{O}_X((m-1)H))+1$ and thus 
\[ \kappa(-K_X)=\kappa(H)=1. \qedhere\]
\end{proof}

\section{Anticanonical system without fixed part}\label{3}
In this section, we consider the following setup: 
\begin{set}\label{setup mobile}
Let $X$ be a smooth projective rationally connected threefold with anticanonical bundle $-K_X$ nef, $n(-K_X)=3$ and $\nu(-K_X)=2$. We suppose that the anticanonical system $|-K_X|$ has no fixed part, so that by Corollary \ref{mH} we can write 
\[
-K_X\sim mH
\]
with $m\geq 2$ and $H$ some prime divisor.
\end{set}

We may now run the minimal model program. Consider an extremal contraction $\varphi\colon X\to Y$.

\subsection{Del Pezzo fibrations}
\begin{prop}\label{delpezzo}
In Setup \ref{setup mobile}, suppose that there exists an extremal contraction $\varphi\colon X\to \mathbb{P}^1$. Then $X\subset \mathbb{P}(\mathcal{E})$ with
\[
\mathcal{E}=\mathcal{O}_{\mathbb{P}^1}\oplus\mathcal{O}_{\mathbb{P}^1}\oplus\mathcal{O}_{\mathbb{P}^1}(-1)\oplus\mathcal{O}_{\mathbb{P}^1}(-1),
\]
and $X\in |\mathcal{O}_{\mathbb{P}(\mathcal{E})}(2) + 4F|$, where $F$ is a general fiber of $\pi\colon\mathbb{P}(\mathcal{E})\to\mathbb{P}^1$.
\end{prop}
\begin{proof}
We use notation from Setup \ref{setup mobile}. Since $-K_X\sim mH$ with $m\geq 2$, we deduce from the classification of Mori-Mukai \cite[Section 3]{MR715648} that $m=2$ or $3$.

{\em Case 1.}
 If $m=3$, then $\varphi$ is a $\mathbb{P}^2$-bundle and we can write $X=\mathbb{P}(\mathcal{E})$ where $\mathcal{E}$ is a vector bundle over $\mathbb{P}^1$ of rank $3$. Denote the tautological line bundle by $\xi\coloneqq\mathcal{O}_{\mathbb{P}(\mathcal{E})}(1)$. Then the Grothendieck relation (see \cite[Appendix A.3]{MR0463157}) gives
$$
\xi^3-\xi^2\cdot\varphi^*(c_1(\mathcal{E}))=0.
$$
Hence
\begin{align*}
    (-K_X)^3 &=(3\xi+\varphi^*(-K_{\mathbb{P}^1}-c_1(\mathcal{E})))^3\\
    &=27\xi^3 + 27\xi^2\cdot\varphi^*(\mathcal{O}_{\mathbb{P}^1}(2)-c_1(\mathcal{E}))\\
    &= 27\xi^2\cdot\varphi^*(\mathcal{O}_{\mathbb{P}^1}(2))\\
    &=54
\end{align*}
which contradicts the fact that $K_X^3=0.$

{\em Case 2.} If $m=2$, then $\varphi\colon X\to \mathbb{P}^1$ is a quadric bundle with general fiber $F_X\simeq \mathbb{P}^1\times\mathbb{P}^1$, and every fiber is a smooth quadric or a quadric cone in $\mathbb{P}^3$. Define
$\mathcal{E}\coloneqq\varphi_*(\mathcal{O}_X(H))$ which is a vector bundle on $\mathbb{P}^1$ of rank $$r=h^0(F_X,H|_{F_X})=h^0(\mathbb{P}^1\times\mathbb{P}^1,\mathcal{O}_{\mathbb{P}^1\times\mathbb{P}^1}(1,1))=4.$$
Now the morphism $\varphi^*\mathcal{E}\to\mathcal{O}_X(H)$ is surjective as it is the evaluation map on each fiber and the restriction of $H$ on each fiber is base-point-free. Hence it gives an embedding $X\subset\mathbb{P}(\mathcal{E})$ such that $H=\mathcal{O}_{\mathbb{P}(\mathcal{E})}(1)|_X$. Let $\pi\colon\mathbb{P}(\mathcal{E})\to \mathbb{P}^1$ such that $\varphi=\pi|_X$.

We write $\mathcal{E}=\oplus_{i=1}^4\mathcal{O}_{\mathbb{P}^1}(a_i)$ with $a_1\geq a_2\geq a_3\geq a_4$. Denote the tautological line bundle by $\xi\coloneqq\mathcal{O}_{\mathbb{P}(\mathcal{E})}(1)$ and the general fiber of $\pi$ by $F$. Since 
\[
K_{\mathbb{P}(\mathcal{E})}= -4\xi + \pi^*(K_{\mathbb{P}^1}+ c_1(\mathcal{E}))
\]
and $K_X=-2\xi|_X$, we deduce from the adjunction formula that
\[
X\in |2\xi + \alpha F|
\]
with $\alpha=-c_1(\mathcal{E})+2,$ because the morphism $\Pic(\mathbb{P}(\mathcal{E}))\to \Pic(X)$ is injective (indeed $\Pic(\mathbb{P}(\mathcal{E}))\simeq \mathbb{Z}\oplus \mathbb{Z}$ and both $\xi$ and $F$ are non-trivial and linearly independent on $X$).

On the other hand, by the Grothendieck relation, we have \[\xi^4-\xi^3\cdot\pi^*(c_1(\mathcal{E}))=0.\] Hence
\[
0=H^3=(\xi|_X)^3=\xi^3\cdot (2\xi+\alpha F)=2c_1(\mathcal{E})+\alpha.
\]
Therefore, $c_1(\mathcal{E})=-2$ and $\alpha=4$.

Since $h^0(\mathbb{P}^1,\mathcal{E})=h^0(X,\mathcal{O}_X(H))=2$ by Corollary \ref{mH}, there are the two following possibilities: either
\[
a_1=1, a_2=a_3=a_4=-1,
\]
or \[a_1=a_2=0, a_3=a_4=-1.\]

Now suppose that $\mathcal{E}=\mathcal{O}(1)\oplus\mathcal{O}(-1)^{\oplus3}$. Then $\Bs|\xi|=\mathbb{P}(\mathcal{O}(-1)^{\oplus3})\eqqcolon D_0$ and $\xi=D_0+F$.

Since $H^0(\mathbb{P}(\mathcal{E}),\mathcal{O}_{\mathbb{P}(\mathcal{E})}(\xi-X))=H^0(\mathbb{P}(\mathcal{E}),\mathcal{O}_{\mathbb{P}(\mathcal{E})}(-\xi-4F))=0$, we deduce from the short exact sequence
\[
0\to \mathcal{O}_{\mathbb{P}(\mathcal{E})}(\xi-X)\to \mathcal{O}_{\mathbb{P}(\mathcal{E})}(\xi)\to \mathcal{O}_X(\xi|_X)\to 0
\]
that the restriction morphism $H^0(\mathbb{P}(\mathcal{E}),\mathcal{O}_{\mathbb{P}(\mathcal{E})}(\xi))\to H^0(X,\mathcal{O}_X(H))$ is injective, hence surjective as $h^0(\mathbb{P}(\mathcal{E}),\mathcal{O}_{\mathbb{P}(\mathcal{E})}(\xi))=h^0(X,\mathcal{O}_X(H))$.

Therefore, when we restrict the base locus $D_0$ of $|\xi|$ to $X$, we have \[D_0\cap X\subset \Bs|H|.\]
But this implies that the base locus of $|H|$ on $X$ has a divisorial part, which contradicts the fact that $|H|$ is mobile on $X$.
\end{proof}

\begin{rem}\label{remdP}
In the setting of Proposition \ref{delpezzo}, $\varphi\colon X\to\mathbb{P}^1$ is a quadric bundle with general fiber $F_X\simeq \mathbb{P}^1\times\mathbb{P}^1$ and $|-K_X|=|2H|$.
Let $D$ be a general member of $|H|$, then $\mathcal{O}_{F_X}(D|_{F_X})\simeq\mathcal{O}_{\mathbb{P}^1\times\mathbb{P}^1}(1,1)$. Hence a general fiber of $\varphi\colon D\to\mathbb{P}^1$ is either isomorphic to $\mathbb{P}^1$ or two $\mathbb{P}^1$'s intersecting transversally at one point.
\end{rem}

\begin{proof}[Proof of Theorem \ref{mainmobile} (1').]
Let $\mathcal{E}=\mathcal{O}_{\mathbb{P}^1}^{\oplus 2}\oplus\mathcal{O}_{\mathbb{P}^1}(-1)^{\oplus 2}$ and $\pi\colon\mathbb{P}({\mathcal{E}})\to\mathbb{P}^1$ be the projection morphism. Denote the tautological line bundle $\Oe(1)$ by $\xi$ and a general fiber of $\pi$ by $F$. Let $X$ be a very general member in $|2\xi+4F|$. Since $\Sym^2\mathcal{E}\otimes\Op(4)$ is globally generated, a general member in $|2\xi+4F|$ is smooth. As
    \[\
    K_{\mathbb{P}(\mathcal{E})}=-4\xi+\pi^*(K_{\mathbb{P}^1}+c_1(\mathcal{E})),
    \]
    the adjunction formula gives
    $K_X=-2\xi|_X.$
    Let $H\coloneqq\xi|_X$, then $-K_X=2H$.
    
    Let $\mathcal{E}_0\coloneqq\Op\oplus\Op(-1)^{\oplus 2}$ and $D_0\coloneqq\mathbb{P}(\mathcal{E}_0)$. Then $D_0\in|\xi|$ and we have the projection morphism $\pi_0\coloneqq\pi|_{D_0}\colon D_0\to\mathbb{P}^1$ and the tautological line bundle $\xi_0=\xi|_{D_0}$ associated to $\mathcal{O}_{\mathbb{P}(\mathcal{E}_0)}(1)$. Let $H_0\coloneqq X\cap D_0$. Then $H_0\in |H|$. Let $S_0\coloneqq\mathbb{P}(\Op(-1)^{\oplus 2})\simeq \mathbb{P}^1\times\mathbb{P}^1$. Then $S_0\in|\xi_0|$.
    
    Since $R^i \pi_*(\Oe(1))=0$ for all $i>0$,
    \begin{align*}
    H^1(\mathbb{P}(\mathcal{E}),\Oe(X-D_0))&\simeq H^1(\mathbb{P}^1,\mathcal{E}\otimes\Op(4))\\
    &= H^1(\mathbb{P}^1, \Op(4)^{\oplus 2}\oplus \Op(3)^{\oplus 2})\\
    &= 0.
    \end{align*}
    Similarly we have 
    \begin{align*}
    H^1(D_0,\mathcal{O}_{D_0}(H_0-S_0))&\simeq H^1(\mathbb{P}^1,\mathcal{E}_0\otimes\Op(4))\\
    &= H^1(\mathbb{P}^1,\Op(4)\oplus \Op(3)^{\oplus 2})\\
    &= 0.
    \end{align*}
    Therefore, the evaluation maps
    \[\ H^0(\mathbb{P}(\mathcal{E}),\Oe(X))\to H^0(D_0,\mathcal{O}_{D_0}(X|_{D_0}))\]
    and
    \[\ H^0(D_0,\mathcal{O}_{D_0}(H_0))\to H^0(S_0,\mathcal{O}_{S_0}(H_0|_{S_0})) \]
    are surjective. 
    On the other hand, since $2\xi+4F$ is globally generated, its restriction to $D_0$ (resp. to $S_0$) is globally generated. Hence by the surjectivity of the above evaluation maps, we deduce that $H_0=X\cap D_0$ (resp. $C_0\coloneqq X\cap S_0$) is smooth for a general $X\in |2\xi+4F|$.
    
    {\em Claim.} $H$ (and thus $-K_X$) is nef.
    
    For any curve $C\subset \mathbb{P}(\mathcal{E})$ such that $\xi\cdot C<0$, we have $C\subset S_0$.
    
    Denote the two ruling of  $S_0$ by $f_1$ and $f_2$, where $f_1\coloneqq F|_{S_0}$ and $f_2$ surjects to $\mathbb{P}^1$ by $\pi$. Then 
    \[\ \xi|_{S_0}=\mathcal{O}_{\mathbb{P}(\Op(-1)^{\oplus 2})}(1)\sim -f_1+f_2. \]
    
    Therefore \[\ X|_{S_0}\sim (2\xi+4F)|_{S_0}\sim 2(-f_1+f_2)+4f_1 = 2(f_1+f_2).\]
   
   Now suppose by contradiction that there exists an integral curve $C\subset X$ such that $-K_X\cdot C<0$. Then $\xi\cdot C=\xi|_X\cdot C<0$ and thus $C\subset X\cap S_0$. But $C_0=X\cap S_0$ is a smooth irreducible curve (it is a smooth elliptic curve), we deduce that \[\ [C]=[2(f_1+f_2)],\] which implies that
   \[\ \xi\cdot C = \xi|_{S_0}\cdot C=(-f_1+f_2)\cdot 2(f_1+f_2)=0. \]
   This contradicts the fact that $\xi\cdot C<0$. Hence $H$ is nef and this proves the claim.
   
    By the adjunction formula, one has that $-K_{H_0}\sim H|_{H_0}$ is nef with $(-K_{H_0})^2=0$ and $C_0\in |-K_{H_0}|$.  Furthermore, $\pi$ induces a fibration on $H_0$ over $\mathbb{P}^1$ with general fiber isomorphic to $\mathbb{P}^1$.
    
    Now let $S$ be the blow-up of $\mathbb{P}^2$ at $9$ points in very general position such that $-K_S$ is nef, not semi-ample and the unique member $D\in|-K_S|$ is a smooth elliptic curve. Denote the blow-up by $\sigma:S\to\mathbb{P}^2$. Let $h=\sigma^*(\mathcal{O}_{\mathbb{P}^2}(1))$ and $C_i$ be a conic on $S$, i.e. a smooth rational curve such that $-K_S\cdot C_i=2$ and $C_i^2=0$ (for example take $C_i$ the strict transform of a general line through a blown-up point $p_i\in\mathbb{P}^2$ such that $C_i\sim h-e_i$ where $e_i$ is the exceptional curve over $p_i$). Then the class of $C_i$ induces a conic bundle $\tau:S\to\mathbb{P}^1$.
    
    Since $\tau\colon S\to\mathbb{P}^1$ is a regular conic bundle, one has
    \[ R^i \tau_*(\mathcal{O}_S(-K_S))=0\]
    for all $i>0$
    and $\tau_*(\mathcal{O}_S(-K_S))$ is a locally free sheaf of rank $3$ that we denote by $\mathcal{V}$. Therefore, 
    \[ H^k(\mathbb{P}^1,\mathcal{V})\simeq H^k(S,\mathcal{O}_S(-K_S))\]
    for all $k\geq 0$ and thus $\chi(\mathbb{P}^1,\mathcal{V})=\chi(S,\mathcal{O}_S(-K_S))=1$.
    Now by the Grothendieck-Riemann-Roch theorem, one has 
    \[ \chi(\mathbb{P}^1,\mathcal{V})=\dg(\mathcal{V})+3\]
   and thus $c_1(\mathcal{V})=-2$.
   
   Now since $h^0(\mathbb{P}^1,\mathcal{V})=h^0(S,\mathcal{O}_S(-K_S))=1$, we can write
   \[\ \mathcal{V}\simeq \Op\oplus \Op(a)\oplus \Op(b)\]
   with $a,b<0$. As $a+b=-2$, we deduce that $a=b=-1$. Therefore
   \[\ S\subset \mathbb{P}(\Op\oplus\Op(-1)^{\oplus 2})=D_0 \]
   and we have $S\in|(2\xi+4F)|_{D_0}|$ by the adjunction formula. Hence by semicontinuity of cohomology, the surface $H_0=X\cap D_0$ has nef and not semi-ample anticanonical divisor for a very general element $X\in|2\xi+4F|$. Thus $-K_X$ is not semi-ample by Lemma \ref{lemconverse}.
   
\end{proof}

\subsection{Conic bundles}
\begin{prop}\label{conicbundle}
In Setup \ref{setup mobile}, suppose that there exists an extremal contraction $\varphi\colon X\to Y$ to a surface $Y$. Then $X=\mathbb{P}(\mathcal{E})$ is a $\mathbb{P}^1$-bundle over $Y$ with $-K_Y$ nef, $\mathcal{E}$ is a nef rank $2$ vector bundle with $c_1(\mathcal{E})=-K_Y$ and $c_2(\mathcal{E})=K_Y^2$, given by an extension
\[
0\to\mathcal{O}_Y\to \mathcal{E}\to\mathcal{I}_Z\otimes \mathcal{O}_Y(-K_Y)\to 0,
\]
where $\mathcal{I}_Z$ is the ideal sheaf of $c_2(\mathcal{E})$ points on $Y$.
Furthermore, the set of such $X$ forms a bounded family.
\end{prop}
\begin{proof}
We use notation from Setup \ref{setup mobile}. By the classification of Mori-Mukai \cite[Section 3]{MR715648}, $\varphi$ is a conic bundle and $Y$ is a smooth rational surface. Since $-K_X\sim mH$ with $m\geq 2$, we deduce from the classification that $m=2$ and $\varphi$ is a $\mathbb{P}^1$-bundle. By \cite[Proposition 3.1]{MR1255695}, the anticanonical bundle $-K_Y$ is nef. Let $d\coloneqq (-K_Y)^2$, we have thus $0\leq d\leq 9$ and $Y$ is isomorphic to $\mathbb{P}^1\times\mathbb{P}^1$, or $\mathbb{F}_2$, or $\mathbb{P}^2$ blown up in $(9-d)$ points.

We write $X=\mathbb{P}(\mathcal{E})$ with $\mathcal{E}=\varphi_*(\mathcal{O}_X(H))$. Then $H=\xi:=\mathcal{O}_{\mathbb{P}(\mathcal{E})}(1)$. As $-K_X=2H$ and 
\[
-K_X=\varphi^*(-K_Y-\det(\mathcal{E}))+2\xi,
\]
one has $c_1(\mathcal{E})=\det(\mathcal{E})=-K_Y$.

On the other hand, since $(-K_X)^3=0$, one has 
\[
0=\xi^3=c_1^2(\mathcal{E})-c_2(\mathcal{E}),
\]
from which we deduce $c_2(\mathcal{E})=K_Y^2=d$.

{\em Claim.} $\mathcal{E}$ has a section which vanishes in codimension at least $2$.

Suppose by contradiction that every non-zero section in $H^0(Y,\mathcal{E})$ vanishes in codimension $1$. Let $s\in H^0(Y,\mathcal{E})$ be a non-zero section and $H_s$ the element in $|H|$ associated to $s$. Let $D$ be the one-dimensional components of the vanishing locus of $s$ taken with multiplicity. Now consider the vector bundle $\mathcal{E}'\coloneqq\mathcal{E}\otimes\mathcal{O}_Y(-D)$. Then it has a non-zero section $s'\in H^0(Y,\mathcal{E}')$ which vanishes in codimension at least $2$. We denote the element associated to $s'$ in $|\mathcal{O}_{\mathbb{P}(\mathcal{E}')}(1)|$ by $H_{s'}$. Then one has an isomorphism $X\simeq \mathbb{P}(\mathcal{E}')$ under which $H_{s'}$ corresponds to $H_s\otimes\varphi^*(-D)$. Hence there exists an effective divisor $R$ on $X$ (which corresponds to $H_{s'}$) such that 
\[
H_s=\varphi^*(D)+ R.
\]
Notice that $R$ is non-zero as the restriction of $H_s$ to a general fiber is $\mathcal{O}_{\mathbb{P}^1}(1)$.
Since this holds for every non-zero section $s\in H^0(Y,\mathcal{E}),$ it contradicts the fact that $H$ is irreducible and reduced.
This proves the claim.

Therefore, following \cite[Section 4.1, page 85--87]{MR1439504}, we have an exact sequence
\[
(*) \quad \quad 0\to\mathcal{O}_Y\to \mathcal{E}\to\mathcal{I}_Z\otimes \mathcal{O}_Y(-K_Y)\to 0,
\]
where $Z$ is the zero locus of a general section of $\mathcal{E}$ with length $l(Z)=c_2(\mathcal{E})=d$.

If $d=0$, then we have $Z=\emptyset$ and $(*)$ must split as \[\Ext^1(\mathcal{O}_Y(-K_Y),\mathcal{O}_Y)\simeq H^1(Y,\mathcal{O}_Y(K_Y))\simeq H^1(Y,\mathcal{O}_Y)=0,\] thus $\mathcal{E}=\mathcal{O}_Y\oplus\mathcal{O}_Y(-K_Y)$. Consider the case when $d>0$. For a fixed smooth rational surface $Y$ such that $-K_Y$ is nef, $Z$ is a finite subscheme of length $d=K_Y^2$ on $Y$. Hence it is parameterized by the Hilbert scheme $Y^{[d]}$. Furthermore, the extensions $(*)$ are parameterized by the vector space $\Ext^1(\mathcal{I}_Z\otimes\mathcal{O}_Y(-K_Y),\mathcal{O}_Y)$ of finite dimension. Therefore, the set of varieties $\mathbb{P}(\mathcal{E})$ such that $\mathcal{E}$ is a vector bundle of rank $2$ over $Y$ satisfying $(*)$ forms a bounded family. 

Now since the set of smooth rational surfaces $Y$ with $-K_Y$ nef forms a bounded family (see \cite[Section 6]{MR1298994}), we deduce that the set of such $X=\mathbb{P}(\mathcal{E})$ forms a bounded family as well.
\end{proof}

\begin{rem}\label{remconic}
In the setting of Proposition \ref{conicbundle}, one has $X=\mathbb{P}(\mathcal{E})$ where $\mathcal{E}$ is a rank-two vector bundle on the surface $Y$ and $-K_X=2H$, where $H$ is the tautological line bundle $\mathcal{O}_{\mathbb{P}(\mathcal{E})}(1)$.

Let $D$ be a general member in $|H|$. Since $\mathcal{E}$ is given by the short exact sequence $(*)$, one has
\[
D=\Bl_Z(Y).
\]
\end{rem}

\begin{ex}\label{example}
\rm{Let $S$ be $\mathbb{P}^2$ blown up in $9$ points in sufficiently general position such that $-K_S$ is nef and not semi-ample. Then there exists a unique element $D\in|-K_S|$. We have $\kappa(-K_S)=0$ and $K_S^2=0$.

Now define $\mathcal{E}\coloneqq \mathcal{O}_S\oplus\mathcal{O}_S(-K_S)$ and $\pi\colon X\coloneqq\mathbb{P}(\mathcal{E})\to S$. Thus $\mathcal{E}$ is nef and $-K_X=2\xi$, where $\xi\coloneqq\mathcal{O}_{\mathbb{P}(\mathcal{E})}(1)$, is nef. Furthermore, we have $c_1(\mathcal{E})=D$ and $c_2(\mathcal{E})=0$.

For $n\in \mathbb{N}^*$, we have 
\[
h^0(X,\mathcal{O}_X(-nK_X))=h^0(S,\Sym^{2n}(\mathcal{E}))=2n+1.
\]
Hence $\kappa(-K_X)=1$.

Now we consider the sections associated to $\pi$. 
Notice that for any extension $0\to L\to \mathcal{E}\to Q\to 0$ where $L$ and $Q$ are line bundles on $S$, we have \[\mathbb{P}(Q)=\xi-\pi^*(L).\] Hence
there are two types of sections: either it corresponds to the quotient $\mathcal{E}\to\mathcal{O}_S(-K_S)\to 0$ and thus gives an element $D_1\simeq S$ such that $D_1\in|\xi|$, or it corresponds to the quotient $\mathcal{E}\to\mathcal{O}_S\to 0$ and thus gives an element $D_2\simeq S$ such that $D_2\in|\xi-\pi^*D|$. Therefore, there are two types of elements in $|\xi|$: one of the form $D_1$ and the other of the form $D_2+\pi^*D$, where $D_1$ and $D_2$ are two disjoint sections of $\pi$.

Since $D_1\in|\xi|$ moves, $D_1^2$ is an effective $1$-cycle. By the Grothendieck relation, one has $\xi^2-\xi\cdot\pi^*c_1(\mathcal{E})=0$. Hence $\xi^2=D_1^2=D_1\cdot \pi^*D$ is a non-zero effective $1$-cycle isomorphic to $D$. Furthermore, $$\xi^3=\xi\cdot(\pi^*D)^2=0,$$ as $D^2=0$. Therefore, $\nu(-K_X)=2$.}
\end{ex}  

\begin{proof}[Proof of Theorem \ref{mainmobile} (2').]
Let $\mathcal{E}$ be a nef rank-two vector bundle on a smooth rational surface $Y$ with nef anticanonical divisor such that $c_1(\mathcal{E})=-K_Y$, $c_2(\mathcal{E})=(-K_Y)^2$, fitting into a sequence
   \[\ 
   0\to\mathcal{O}_Y\to \mathcal{E}\to\mathcal{I}_Z\otimes \mathcal{O}_Y(-K_Y)\to 0\quad (*)
   \]
where $\mathcal{I}_Z$ is the ideal sheaf of $c_2(\mathcal{E})$ points in sufficiently general position.
Let $X=\mathbb{P}(\mathcal{E})$ and $\xi$ be the tautological line bundle. Let $H$ be a general member in $|\xi|$. Then 
\[\ -K_X\sim 2H\]
and $(-K_X)^3=8\xi^3=8(c_1(\mathcal{E})^2-c_2(\mathcal{E}))=0$.
Furthermore,
the sequence $(*)$ gives 
\[\ H\simeq \Bl_Z(Y),\quad N_{H/X}=-K_H.\]

Since $Y$ is a smooth rational surface with $-K_Y$ nef, $Y$ is isomorphic to $\mathbb{P}^1\times\mathbb{P}^1$, or $\mathbb{F}_2$, or the blow-up of $\mathbb{P}^2$ in $9-(-K_Y)^2$ points in almost general position. 
Note that the blow-up of $\mathbb{P}^1\times\mathbb{P}^1$ or $\mathbb{F}_2$ in a general point is isomorphic to $\mathbb{P}^2$ blown up in $2$ points (see \cite[page 13]{MR2427467}).
Since $Z\subset Y$ is the subscheme of $(-K_Y)^2$ points in sufficiently general position, $H$ is isomorphic to $\mathbb{P}^2$ blown up in at most $9$ points in sufficiently general position. Therefore, $-K_H$ is nef and not semi-ample. Hence $-K_X$ is nef and not semi-ample by Lemma \ref{lemconverse}.
\end{proof}

\subsection{Birational contractions}
\begin{prop}\label{birat}
In Setup \ref{setup mobile}, suppose that there exists a birational extremal contraction $\varphi\colon X\to Y$. Then $-K_X$ is divisible by $2$ in $\Pic(X)$, $Y$ is a smooth almost del Pezzo threefold of degree $1$ and $\varphi$ is the blow-up of a point $p\in Y$. Furthermore, if we write $-K_Y\sim 2H_Y$, then $p\not\in \Bs|H_Y|$.
\end{prop}
\begin{proof}
We use notation from Setup \ref{setup mobile}. Since $-K_X\sim mH$ with $m\geq 2$, by the classification of Mori-Mukai contractions on smooth threefolds \cite[Section 3]{MR715648}, one has $m=2$ and $\varphi$ is the blow-up of a smooth point $p$ on $Y$ with exceptional divisor $E\simeq\mathbb{P}^2$ and $\mathcal{O}_E(E)=\mathcal{O}_{\mathbb{P}^2}(-1)$. Hence $-K_Y$ is nef by \cite[Proposition 3.3]{MR1255695} and 
\[
(-K_Y)^3 = (-K_X)^3 + (2E)^3 = 8,
\]
i.e. $-K_Y$ is big.

On the other hand,
\[
-K_Y = \varphi_*(-K_X) = 2\varphi_*(H)\eqqcolon 2H_Y
\]
with $H_Y\in \Pic(Y)$. Then $H_Y$ is nef and big with $(H_Y)^3=1$. We conclude that $Y$ is an almost del Pezzo threefold of degree $1$ and the base scheme of $|H_Y|$ is one point by \cite[Section 2]{MR2427467}.

If $p$ is the base point of $|H_Y|$, then $\Bs|H|=\emptyset$ since the base scheme of $|H_Y|$ is one point. This is absurd because $|-K_X|$ is not semi-ample.
\end{proof}

\begin{rem}\label{rembirat}
In the setting of Proposition \ref{birat}, let $D\in |H|$ be a general member and $E$ the exceptional divisor of $\varphi$, one has
\[
K_E=(K_X+E)|_E=(-2D+E)|_E
\]
by the adjunction formula.

Since $D$ moves, $D\cdot E$ is an effective $1$-cycle. We deduce from $N_{E/X}\simeq \mathcal{O}_{\mathbb{P}^2}(-1)$ that $D\cdot E=l$, where $l$ is a line on $E\simeq \mathbb{P}^2$.

On the surface $D$, one has
\[
l^2=(E|_D)^2=D\cdot E^2=D|_E\cdot E|_E=-1
\]
and 
\[
K_D\cdot l=(K_X+D)|_D\cdot l=-D\cdot l=-1.
\]
Hence $l$ is a $(-1)$-curve on $D$.
\end{rem}

\begin{proof}[Proof of Theorem \ref{mainmobile} (3').]
Let $Y$ be a smooth almost del Pezzo threefold of degree one. Then a general member in $|-\frac{1}{2}K_Y|$ is a smooth almost del Pezzo surface of degree one.
Now fix a general member $H_Y\in|-\frac{1}{2}K_Y|$. Since $H_Y$ is $\mathbb{P}^2$ blown up at $8$ points in almost general position, by choosing a sufficiently general point $p\in H_Y\subset Y$, the blow-up $D\coloneqq \Bl_p H_Y$ of $H_Y$ at $p$ has nef and not semi-ample anticanonical divisor. Then
let $\varphi\colon X\to Y$ be the blow-up at $p$. We have \[\ -K_X=2(\varphi^*H_Y-E)\eqqcolon 2H\]
and $D\in|H|$. Therefore, we deduce by Lemma \ref{lemconverse} that $-K_X$ is nef and not semi-ample.
\end{proof}

\section{Anticanonical system with non-empty fixed part}\label{4}
We consider the case when the anticanonical system $|-K_X|$ has a non-empty fixed part, that is, we can write $|-K_X|=A+|B|$ with $A$ the fixed part and $|B|$ the mobile part. 
By Corollary \ref{mH}, we have $|B|=|mH|$ where $m\geq 2$ and $H$ is some prime divisor.

\begin{prop}\label{nefmobile}
Let $X$ be a smooth projective rationally connected threefold with $-K_X$ nef, $n(-K_X)=3$, $\nu(-K_X)=2$. If the anticanonical system \[|-K_X|=A+|mH|,\quad m\geq 2\] has non-empty fixed part $A$, then
there exists a finite sequence of flops $\psi\colon X\dashrightarrow X'$ such that $X'$ is smooth with $-K_{X'}$ nef and $H'\coloneqq \psi_*(H)$ is nef.
\end{prop}
\begin{proof}
Fix a general member $F\in |H|$. Since $X$ is smooth, for sufficiently small $\epsilon>0$, the pair $(X,\epsilon F)$ is log-canonical. It follows from Lemma \ref{nef} that if $F$ is not nef, then there exists a $(K_X+\epsilon F)$-negative extremal ray $\Gamma$ such that $\epsilon F\cdot \Gamma<0$. Let $c_{\Gamma}$ be the contraction of the extremal ray $\Gamma$ and $l$ a contracted curve. Thus $F\cdot l<0$, which implies $l\subset \Bs(|H|)$. Since $|H|$ is mobile, it follows that $c_{\Gamma}$ is small. This implies that $K_X\cdot l=0$ since there is no flipping contraction for smooth threefolds. Hence there exists a flop of $c_{\Gamma}$ and the flopped threefold $X^+$ is smooth by \cite[Theorem 2.4]{MR986434}. 

By repeating the above argument and the termination of three-dimensional flops \cite[Corollary 6.19]{MR1658959}, we deduce that there exists a sequence of flops $\psi\colon X\dashrightarrow X'$ such that $H'\coloneqq \psi_*(H)$ is nef.
\end{proof}

\begin{lem}
In the setting of Proposition \ref{nefmobile}, if $H$ is nef, then \[A^3=A^2\cdot H=A\cdot H^2=H^3=0.\]
\end{lem}
\begin{proof}
As $-K_X\sim A+mH$ is nef, one has $K_X^2\cdot A\geq 0$ and $K_X^2\cdot H\geq 0$. Then
\[
0=(-K_X)^3=K_X^2\cdot (A+mH)
\]
gives $K_X^2\cdot A=K_X^2\cdot H=0$. From this we further conclude that
\[
0=-K_X\cdot (A+mH)\cdot H=-K_X\cdot (A\cdot H+mH^2).
\]
Since $H$ moves, $A\cdot H$ and $H^2$ are effective cycles. This implies that \[-K_X\cdot A\cdot H=-K_X\cdot H^2=0.\]
Hence, $A^2\cdot H+mA\cdot H^2=0$ and $A\cdot H^2+ mH^3=0$. As $H$ is nef, $A\cdot H$ and $H^2$ are effective cycles, we deduce that \[A\cdot H^2=H^3=0.\] This implies $A^2\cdot H=0.$ Together with $K_X^2\cdot A=0$, we conclude that $A^3=0.$
\end{proof}

After performing possibly a sequence of flops, the mobile part $|B|=|mH|$ of the anticanonical system $|-K_X|$ becomes nef. In this case, either $B^2=0$ and we are in the case described in \cite[Proposition 7.2]{MR2129540}, or $B^2$ is a non-zero effective $1$-cycle which is the case we will study in the rest of the section.

\subsection{Description of the anticanonical system}
\begin{prop}\label{F}
Consider as above a smooth projective rationally connected threefold $X$ with $-K_X$ nef, $n(-K_X)=3$, $\nu(-K_X)=2$. Suppose that the anticanonical system $|-K_X|=A+|mH|$, $m\geq 2$ has non-empty fixed part $A$, and $H$ is nef such that $H^2$ is a non-zero effective $1$-cycle. Let $F$ be a general member of $|H|$. Then $-K_F$ is nef, effective and divisible by $r\geq 2$ in $\NS(F)$. Furthermore, $\kappa(F,-K_F)=0$, $K_F^2=0$ and $F$ is not covered by $(-K_F)$-trivial curves.
\end{prop}
\begin{proof} By the adjunction formula, we have
\[
-K_F=-(K_X+F)|_F = A|_F + (m-1)F|_F. 
\]
As $F$ is nef, it suffices to show that $A|_F$ is nef:
suppose that there exists an irreducible curve $l\subset F$ such that $A|_F\cdot l<0$. Then $l$ is an irreducible component of the effective cycle $C\coloneqq A|_F$. On the other hand, $F$ is nef and $F\cdot C=0$ as $A\cdot F^2=0$, from which we deduce that $F\cdot l=0.$ Hence
\[
-K_X\cdot l=A\cdot l +  mF\cdot l = A\cdot l <0,
\]
which contradicts the fact that $-K_X$ is nef. Therefore, the restriction $A|_F$ is nef.

We note that $A|_F$ cannot be zero: since $-K_X$ is nef with numerical dimension two, the support of a divisor $D\in |-K_X|$ is connected in codimension one by \cite[Lemma 2.3.9]{MR1668575}.

Now let $\nu\colon\tilde{F}\to F$ be a desingularization of the surface $F$. Since $A|_F$ and $F|_F$ are nef Cartier divisors such that $A|_F\cdot F|_F=0$, their pullbacks to the desingularization $\tilde{F}$ remain nef and orthogonal to each other.
Let \[V\coloneqq <\nu^*(A|_F),\nu^*(F|_F)>\subset \NS(\tilde{F}).\] Let $H$ be an ample divisor on $\tilde{F}$, then $NS(\tilde{F})=\mathbb{R}H\oplus(\mathbb{R}H)^{\perp}$. If ${\rm{dim}}\ V\geq 2$, then ${\rm{dim}}\ (V\cap (\mathbb{R}H)^{\perp})\geq 1$. Hence there exists $v\in V\cap (\mathbb{R}H)^{\perp}$ which is non zero, and $v^2<0$ by the Hodge index theorem. But $v=\lambda\nu^*(A|_F)+\mu\nu^*(F|_F)$ with $\lambda,\mu\in\mathbb{R}$, which implies $v^2\geq 0$. This is absurd. Hence ${\rm{dim}}\ V=1$, i.e. $\nu^*(A|_F)$ and $\nu^*(F|_F)$ are non-zero and numerically proportional. Hence $-K_F$ is divisible by $r\geq 2$ with $r\in \mathbb{N}$.

The surface $F$ is not covered by $(-K_F)$-trivial curves: otherwise, $F$ is covered by $(-K_X)$-trivial curves as $-K_F=-K_X|_F - F|_F$ and $-K_F$ is numerically proportional to $F|_F$. As $F$ moves in $X$, this implies that $X$ is covered by $(-K_X)$-trivial curves. This is absurd because $n(-K_X)=3$.

Furthermore, as $A^2\cdot F=A\cdot F^2=F^3=0$, we have $K_F^2=0$. 

It remains to show that $\kappa(F,-K_F)=0$. Indeed, for any $n\in\mathbb{N}$, we have
\[
1\leq h^0(F, \mathcal{O}_F(-nK_F))\leq h^0(\tilde{F},\nu^*(-nK_F)).
\]
If $h^0(\tilde{F},\nu^*(-nK_F))\geq 2$ for some $n$, then the linear system $|\nu^*(-nK_F)|$ has some non-zero mobile part $M$ on $\tilde{F}$, and $\nu^*(-K_F)\cdot M=0$ as $(-K_F)^2=0$ and $-K_F$ is nef. Hence $\tilde{F}$ is covered by $\nu^*(-K_F)$-trivial curves, from which we deduce that $F$ is covered by $(-K_F)$-trivial curves. This is absurd.
\end{proof}

In order to get a more precise description on the geometric structure of $A$ and $F$, we need the two following lemmas:
\begin{lem}\label{baselemma}
Let $S$ be a projective Gorenstein surface such that the anticanonical divisor $-K_S$ is of the following form:
\[
-K_S=D_1+D_2,
\]
where $D_1$ is effective, $D_2$ is a non-zero effective Cartier divisor which is nef and divisible by $r\geq 2$ in $\NS(S)$.

Suppose that $D_2^2=0$ and that one of the following assertions holds:
\begin{enumerate}[label=\normalfont(\roman*)]
    \item $S$ is not covered by $D_2$-trivial curves;
    \item $D_2$ contains a smooth curve of positive genus.
\end{enumerate}
Then $D_1=0$ and $S$ is normal with at most rational singularities. Furthermore, the surface $\tilde{S}$ obtained by the minimal resolution of $S$ is relatively minimal.
\end{lem}
\begin{proof}
{\em Special case.}
Assume that $S$ is smooth. Suppose by contradiction that $D_1$ is not zero. Since $D_2$ is divisible by $r\geq 2$, we put $rL:\equiv_{num}D_2$, with $L$ nef and $L^2=0$. Then
\[
-(K_S+rL)=D_1
\]
is effective. We deduce that the adjoint bundle $K_S+rL$ is not nef.

Now for every irreducible reduced curve $l\subset S$ such that $(K_S+ rL)\cdot l<0$, one has $K_S\cdot l<0$ since $L$ is nef. Then by the cone theorem, there exists a $K_S$-negative extremal ray $R$ which is $(K_S+ rL)$-negative. We denote the contraction of the extremal ray $R$ by $\phi\colon S\to Z$.

   \begin{enumerate}
        \item If ${\rm{dim}}\ Z=1$, then $\phi\colon S\to Z$ is a $\mathbb{P}^1$-bundle over a smooth algebraic curve $Z$. Let $f$ be a fiber of $\phi$, then $f^2=0$ and $(K_S+rL)\cdot f<0$.
        
        Since $K_S\cdot f=-2$ by the adjunction formula, together with $L\cdot f\geq 0$ and $r\geq 2$, we have $L\cdot f=0$.
        
        \item If $Z$ is a point, then $S=\mathbb{P}^2$. But $L$ is nef, not ample and not numerically trivial, this is absurd.
\end{enumerate}

If $\phi$ is birational, let $l\in R$ be an integral contracted curve, then $l$ is actually a $(-1)$-curve since we contract a $K_S$-negative extremal ray. Hence $L\cdot l=0$ as $(K_S+rL)\cdot l<0$ and $L$ is nef.
    
    Now we put $L'\coloneqq \phi_*(L)$ and $D_1'\coloneqq \phi_*(D_1)$. Then
    \[
    -K_Z=rL'+D_1'.
    \]
    Since $L\cdot l=0$, we know by the contraction theorem that $L\simeq \phi^*(L')$.
    Hence $L'$ is nef and $L'^2=0$.
    
    Notice that the two assertions in the lemma are preserved by the contraction $\phi$. More precisely, 
     \begin{enumerate}[label=(\roman*)]
     \item if $Z$ is covered by $\phi_*(D_2)$-trivial curves, then $S$ is covered by $D_2$-trivial curves as $L=\phi^*(L')$;
     \item $\phi_*(D_2)$ contains a curve of positive genus, as $\phi$ does not contract any curve of positive genus.
     \end{enumerate}
    
    Moreover, $Z$ cannot be a minimal surface. Indeed if $K_Z$ is nef, then \[K_S+rL=\phi^*(K_Z)+C+rL\] is pseudo-effective. Therefore, $K_S+rL$ is zero as it is anti-effective. This is absurd because $D_1$ is not zero.
    
    Therefore, by running a $(K_S+rL)$-minimal model program, we can suppose that $S$ is a $\mathbb{P}^1$-bundle as described in the first case above. Now we show that this will lead to a contradiction:
    \begin{enumerate}[label=(\roman*)]
    \item We first consider the case when $S$ is not covered by $D_2$-trivial curves: since $L\cdot f=0$ for every fiber $f$ of $\phi$, the surface $S$ is covered by $L$-trivial curves. Hence $S$ is covered by $D_2$-trivial curve, which is absurd.
    \item For the case when $D_2$ contains a smooth curve of positive genus: since $$D_2\cdot f=rL\cdot f=0$$ for a general fiber $f$ of $\phi$, $D_2$ is contained in some special fiber of the $\mathbb{P}^1$-bundle. This is absurd because $D_2$ contains a curve of positive genus.
    \end{enumerate}

Therefore, $D_1=0$. Furthermore, since $-K_S=D_2$ is divisible by $r\geq 2$, the surface $S$ does not contain any $(-1)$-curve, i.e. $S$ is relatively minimal.

{\em General case.}
Let $\nu\colon\overline{S}\to S$ be the normalization of $S$ and $\mu\colon\tilde{S}\to \overline{S}$ the minimal resolution of $\overline{S}$. We put $\pi\coloneqq\nu\circ\mu \colon\tilde{S}\to S$. Computing the anticanonical bundles we get
\[
-K_{\overline{S}}=\nu^*(-K_S)+E_1
\]
with some effective Weil divisor $E_1$ supported on the zero locus of the conductor ideal and
\[
-K_{\tilde{S}}=\pi^*(-K_S)+\tilde{E_1}+E_2
\]
with $\tilde{E_1}$ the proper transform of $E_1$ in $\tilde{S}$ and $E_2$ some effective divisor supported on the exceptional locus.

Now $\tilde{S}$ is a smooth surface such that 
\[
-K_{\tilde{S}}=\tilde{D_1}+\tilde{D_2}
\]
with $\tilde{D_1}\coloneqq\tilde{E_1}+E_2+\pi^*(D_1)$ effective divisor, $\tilde{D_2}\coloneqq\pi^*(D_2)$ non-zero, effective, nef and divisible by $r\geq 2$.

Furthermore, one has $\tilde{D_2}^2=D_2^2=0$ and $\tilde{D_2}$ satisfies one of the two assertions in the lemma if $D_2$ does:
\begin{enumerate}[label=(\roman*)]
     \item if $\tilde{S}$ is covered by $\pi^*(D_2)$-trivial curves, then $S$ is covered by $D_2$-trivial curves;
     \item $\pi^*(D_2)$ contains a smooth curve of positive genus which surjects to the one contained in $D_2$.
     \end{enumerate}

Hence by the previous smooth case, we deduce that $\tilde{D_1}=0$.
This implies that $D_1=0$ and $S$ is normal as it is Cohen-Macaulay, with at worst rational singularities.

Let $\mu\colon\tilde{S}\to S$ be the minimal resolution of $S$. Then $-K_{\tilde{S}}=\mu^*(-K_S)=\mu^*(D_2)$ is divisible by $r\geq 2$. Thus $\tilde{S}$ does not contain any $(-1)$-curve, i.e. $\tilde{S}$ is relatively minimal.
\end{proof}  

\begin{lem}(\cite{MR2129540}, Proposition $1.5$, Proposition $1.6$)\label{classification}
Let $S$ be a smooth projective surface with $-K_S$ nef and $\nu(-K_S)=1$. Then $S$ is one of the following:
\begin{enumerate}[label=\normalfont\arabic*.]
    \item $n(-K_S)=1$: $S$ admits an elliptic fibration and $-K_S$ is semi-ample;
    \item $n(-K_S)=2$: we have $\kappa(-K_S)=0$ and either
    \begin{enumerate}[label=\normalfont(\Alph*)]
        \item $S$ is $\mathbb{P}^2$ blown up in $9$ points in sufficiently general position or
        \item $S=\mathbb{P}(\mathcal{E})$ with $\mathcal{E}$ a rank $2$ vector bundle over an elliptic curve which is defined by an extension
$$
0\to \mathcal{O}\to \mathcal{E}\to\mathcal{L}\to 0
$$
with $\mathcal{L}$ a line bundle of degree $0$ and either
\begin{enumerate}[label=\normalfont(\roman*)]
    \item $\mathcal{L}=\mathcal{O}$ and the extension is non-split or
    \item $\mathcal{L}$ is not torsion.
\end{enumerate}
The structure of the unique element $D$ in $|-K_S|$ is as follows:
\begin{enumerate}[label=\normalfont(\roman*)]
    \item $D=2C$ and $C$ is a smooth elliptic curve.
    \item $D=C_1+C_2$ where $C_1$ and $C_2$ are smooth elliptic curves which do not meet. 
\end{enumerate}
    \end{enumerate}
\end{enumerate}
\end{lem}

\begin{cor}\label{corclassification}
In the setting of Lemma \ref{baselemma}, the surface $S$ is smooth. It is a $\mathbb{P}^1$-bundle over a smooth elliptic curve.
\end{cor}
\begin{proof} Let $\mu\colon\tilde{S}\to S$ be the minimal resolution of $S$. Then by Lemma \ref{baselemma} $-K_{\tilde{S}}=\mu^*(-K_S)$ is non-zero, effective and nef. Hence $\tilde{S}$ is uniruled and thus it admits a Mori fibration. Furthermore, since $\tilde{S}$ is relatively minimal by Lemma \ref{baselemma}, we deduce that $\tilde{S}$ is a $\mathbb{P}^1$-bundle over a smooth curve.

Now by the classification in Lemma \ref{classification}, $\tilde{S}$ is either an elliptic fibration or a $\mathbb{P}^1$-bundle over a smooth elliptic curve. In both of the two cases, we deduce that $\tilde{S}$ is a $\mathbb{P}^1$-bundle over a smooth elliptic curve and $S=\tilde{S}$.
\end{proof}

\begin{lem}\label{lemmaF}
Let $X$ be a smooth projective rationally connected threefold $X$ with $-K_X$ nef, $n(-K_X)=3$, $\nu(-K_X)=2$. Suppose that the anticanonical system $|-K_X|=A+|mH|$, $m\geq 2$ has non-empty fixed part $A$, and $H$ is nef such that $H^2$ is a non-zero effective $1$-cycle. Let $F$ be a general member of $|H|$. Then $F$ is a smooth surface such that $-K_F$ is nef and divisible by $2$ in $\NS(F)$ with $\nu(-K_F)=1$, $\kappa(-K_F)=0$. More precisely, $F=\mathbb{P}(\mathcal{E})$ with $\mathcal{E}$ a rank-$2$ vector bundle over an elliptic curve as described in the Lemma \ref{classification}, 2.(B). Furthermore, we have $m=2$ and $A\cdot F$ is a smooth elliptic curve.
\end{lem}
\begin{proof}
By Proposition \ref{F}, we have that $-K_F$ is non-zero, effective, nef and divisible by $r\geq 2$. Furthermore, $(-K_F)^2=0$ and $F$ is not covered by $(-K_F)$-trivial curves. Hence we can apply Lemma \ref{baselemma} and Corollary \ref{corclassification} to obtain that the surface $F$ is a $\mathbb{P}^1$-bundle over a smooth elliptic curve. Now since $F$ is not covered by $(-K_F)$-trivial curves, i.e. $n(-K_F)=2$, we deduce from the classification in Lemma \ref{classification} that $F=\mathbb{P}(\mathcal{E})$ with $\mathcal{E}$ a rank $2$ vector bundle over an elliptic curve defined as in the case $(B)$.

Since $-K_F=A|_F+(m-1)F|_F$, we deduce from the structure of the unique element in $|-K_F|$ that $m=2$ and $A\cdot F$ is a smooth elliptic curve.
 \end{proof}

\begin{lem}\label{lemmaA}
Let $X$ be a smooth projective rationally connected threefold $X$ with $-K_X$ nef, $n(-K_X)=3$, $\nu(-K_X)=2$. Suppose that the anticanonical system $|-K_X|=A+|mH|$, $m\geq 2$ has non-empty fixed part $A$, and $H$ is nef such that $H^2$ is a non-zero effective $1$-cycle. Then $A$ is an irreducible reduced smooth surface such that $-K_A$ is nef and divisible by $2$ in $\Pic(A)$ with $\nu(-K_A)=1$. More precisely, the surface $A$ is a $\mathbb{P}^1$-bundle over a smooth elliptic curve.
\end{lem}
\begin{proof} Let $F$ be a general member in $|H|$. As $A|_F$ is an irreducible reduced curve by Lemma \ref{lemmaF}, we can find a divisor $A_1$ which occurs in $A$ with multiplicity one and the rest $A'$ does not meet $F$. Since $m=2$ and $A\cdot F$ is a smooth elliptic curve by the Lemma \ref{lemmaF}, the adjunction formula gives
\[
-K_{A_1}=(A'+2F)|_{A_1}=A'|_{A_1}+2C_0, 
\]
where $C_0$ is a smooth elliptic curve and $A'|_{A_1}$ is an effective divisor on $A_1$.

Moreover, since $F$ is nef and $A\cdot F^2=0$, $C_0$ is nef and $C_0^2=0$ on $A_1$.

Now we can apply Lemma \ref{baselemma} and Corollary \ref{corclassification} to the surface $A_1$, which gives $A'|_{A_1}=0$ and $A_1$ is a $\mathbb{P}^1$-bundle over a smooth elliptic curve. Moreover, the support of a divisor $D\in |-K_X|$ is connected in codimension one by \cite[Lemma 2.3.9]{MR1668575}. As $A'$ does not meet $F$ and $A'|_{A_1}=0$, we obtain $A'=0$. Thus $A=A_1$ and $-K_A=2F|_A$.
\end{proof}

\begin{proof}[Proof of Proposition \ref{descrip}.] It follows from Lemma \ref{lemmaF} and Lemma \ref{lemmaA}.
\end{proof}

\subsection{Running the minimal model program} 
In this subsection, we consider the following setup: 
\begin{set}\label{setupfixed}
Let $X$ be a smooth projective rationally connected threefold $X$ with $-K_X$ nef, $n(-K_X)=3$, $\nu(-K_X)=2$. Suppose that the anticanonical system $|-K_X|=A+|mH|$, $m\geq 2$ has non-empty fixed part $A$, and $H$ is nef such that $H^2$ is a non-zero effective $1$-cycle.
\end{set}
Remind that in this setup, one has
 $|-K_X|=A+|2H|$, both $A$ and a general member $F$ in $|H|$ are $\mathbb{P}^1$-bundles over a smooth elliptic curve such that their anticanonical divisors are nef and divisible by two in $\Pic(A)$ (resp. in $\NS(F)$). Furthermore, both $A\cdot F$ and $F^2$ are smooth elliptic curves.

Consider an extremal contraction $\varphi\colon X\to Y$.
Let $R$ be the extremal ray contracted by $\varphi$. Recall that the length of an extremal ray $R$ is defined by
\[
l(R)=\min\{-K_X\cdot Z\mid [Z]\in R\}.
\]
Let $l$ be a rational curve such that $[l]\in R$ and $-K_X\cdot l=l(R).$ In the birational case, we denote the exceptional divisor of $\varphi$ by $E$.

\subsubsection{Non-birational cases}
In this part, we will show that the contraction $\varphi\colon X\to Y$ cannot be of Mori fiber type.

\medskip

\noindent
{\em Case ${\rm{dim}}\ Y=1$.}
In this case, $-K_X\cdot l=1,2$ or $3$. Recall that for an extremal contraction $\varphi\colon X\to\mathbb{P}^1$, all the fibers are irreducible. Since $A$ is the fixed part of $|-K_X|$, it cannot be a fiber of $\varphi$. As for $H$, since $H^2$ is a non-zero effective cycle, it cannot be a fiber of $\varphi$. We deduce that $A\cdot l>0$ and $H\cdot l> 0$, as the Picard group of $X$ is generated by a fiber of $\varphi$ and another element which has positive intersection with $l$. Therefore, $-K_X\cdot l=3$, $A\cdot l=H\cdot l=1$ and $\varphi$ is a $\mathbb{P}^2$-bundle over $\mathbb{P}^1$.

Now we can write $X=\mathbb{P}(\mathcal{E})$ with $\mathcal{E}$ a rank-$3$ vector bundle over $Y=\mathbb{P}^1$. After twisting $\mathcal{E}$ by some line bundle, we can suppose that $\mathcal{E}=\varphi_*\mathcal{O}_X(H)$ and $H=\mathcal{O}_{\mathbb{P}(\mathcal{E})}(1)$. Since $H$ is nef, the vector bundle $\mathcal{E}$ is nef. From the fact that a vector bundle on $\mathbb{P}^1$ is nef if and only if it is generated by its global sections, we deduce that $\mathcal{E}$ is generated by its global sections. Therefore, $H=\mathcal{O}_{\mathbb{P}(\mathcal{E})}(1)$ is also generated by its global sections. Since $h^0(X,\mathcal{O}_X(H))=2$ by Corollary \ref{mH} and $H^2\neq 0$, we have $\Bs|H|\neq \emptyset$ which leads to a contradiction.

\medskip

\noindent
{\em Case ${\rm{dim}}\ Y=2$.}
In this case, $\varphi\colon X\to Y$ is a conic bundle and we have $-K_X\cdot l=1$ or $2$.
\begin{enumerate}[label=(\roman*)]
    \item If $F\cdot l=0$, then we have $F=\varphi^*(C)$ where $C$ is an irreducible curve on $Y$. Hence $F^2$ an effective cycle contained in some fiber of $\varphi$. This is absurd because $F^2$ is a smooth elliptic curve.
    \item If $F\cdot l=1$, then $\varphi$ is a $\mathbb{P}^1$-bundle and induces a birational morphism from $F$ to $Y$. This is impossible since $q(F)=1$ and $q(Y)=0$.
\end{enumerate}

\subsubsection{Birational contractions}
Since $X$ is a smooth threefold, the contraction $\varphi$ is divisorial.

\medskip

\noindent
{\em Case $A\cdot l=0$.} 
In this case, we have $F\cdot l=1$ and $-K_X\cdot l=A\cdot l+2F\cdot l=2$. Hence $\varphi$ is the blow-up of a smooth point on $Y$, with exceptional divisor $E\simeq\mathbb{P}^2$ and $N_{E/X}\simeq \mathcal{O}_{\mathbb{P}^2}(-1)$. Now the adjunction formula $K_E=(K_X+E)|_E$ gives
\[
\mathcal{O}_E(A|_E) \otimes\mathcal{O}_E(2F|_E) = \mathcal{O}_E(2).
\]
As $E\cdot F$ is a non-zero effective cycle, we deduce that $A\cdot E=0$ and $E\cdot F=l$. On the other hand, we have
\[
(E|_F)^2=F\cdot E^2= F|_E\cdot (-l)=-1.
\]
Hence $l$ is a $(-1)$-curve on the surface $F$, which contradicts the fact that $F$ is relatively minimal.

\medskip

\noindent
{\em Case $A\cdot l<0$.} 
Since the contraction is divisorial, we have $E=A$ in this case. Since $A$ is a ruled surface over a smooth elliptic curve, we know that $l$ is a fiber of $A$ and $F\cdot l=1$. Therefore, $\varphi$ is the blow-up of an elliptic curve and $Y$ is smooth with $-K_Y$ nef by \cite[Theorem 3.8]{MR1255695}.
Furthermore, as we contract the curves meeting $F$ transversally, we conclude that $G\coloneqq \varphi(F)\simeq F$. Since \[-K_Y=\varphi_*(-K_X)=\varphi_*(A+2F)=2\varphi_*(F)=2G,\]
we see that $|-K_Y|=|2G|$ is without fixed part.

We can compute the Kodaira dimension and the numerical dimension for $-K_Y$:
\[
\kappa(-K_Y)=\kappa(\varphi^*(-K_Y))=\kappa(-K_X + E),
\]
and similarly for the numerical dimension we have
\[
\nu(-K_Y)=\nu(-K_X + E).
\]
On the other hand, since $E=A,$ we have
\[
\kappa(-K_X)\leq \kappa(-K_X+E)\leq \kappa(-2K_X)=\kappa(-K_X)
\]
and similarly
\[
\nu(-K_X)\leq \nu(-K_X+E)\leq \nu(-2K_X)=\nu(-K_X),
\]
we deduce that $\kappa(-K_Y)=\kappa(-K_X)=1$ and $\nu(-K_Y)=\nu(-K_X)=2$.

\medskip

\noindent
{\em Case $A\cdot l>0$.} 
In this case, $F\cdot l=0$ since otherwise $-K_X\cdot l>2$ which is in contradiction to \cite[Section 3]{MR715648}. Furthermore, $E\neq A$ and thus $A\cdot E$ is an effective cycle. We will show that the only possible case is when $\varphi$ contract $E$ to a smooth curve of positive genus.

By the classification of Mori-Mukai \cite[Section 3]{MR715648}, we are in one of the following cases:

\begin{enumerate}[label=(\arabic*)]
\item If $E$ is contracted to a point, then one of following cases occurs:
\begin{enumerate}[label=(\roman*)]
\item $E\simeq \mathbb{P}^2$, $N_{E/X}\simeq\mathcal{O}_{\mathbb{P}^2}(-1)$. In this case, we have $A\cdot l=2$ and the adjunction formula $K_E=(K_X+E)|_E$ gives
\[
\mathcal{O}_E(A|_E) \otimes\mathcal{O}_E(2F|_E) = \mathcal{O}_E(2).
\]
As $A\cdot E$ is a non-zero effective cycle, we deduce that $F|_E=0$ and $A|_E=\mathcal{O}_E(2)$.

\item $E\simeq \mathbb{P}^1\times\mathbb{P}^1$, $N_{E/X}\simeq \mathcal{O}_{\mathbb{P}^1\times\mathbb{P}^1}(-1,-1)$. In this case, we have $A\cdot l=1$ and the adjunction formula gives
\[
\mathcal{O}_E(A|_E) \otimes\mathcal{O}_E(2F|_E) =\mathcal{O}_E(1,1).
\]
As $A\cdot E$ is a non-zero effective cycle, we deduce that $F|_E=0$ and $A|_E=\mathcal{O}_E(1,1)$.

\item $E$ is a quadric cone in $\mathbb{P}^3$ with $N_{E/X}\simeq \mathcal{O}_{\mathbb{P}^3}(-1)\otimes \mathcal{O}_E$. In this case, we have $A\cdot l=1$ and the adjunction formula gives
\[
\mathcal{O}_E(A|_E) \otimes\mathcal{O}_E(2F|_E) =\mathcal{O}_{\mathbb{P}^3}(1)\otimes\mathcal{O}_E=\mathcal{O}_E(2l).
\]
But since $F|_E$ is Cartier, one cannot have $F|_E=l$ which is $2$-Cartier.
Hence $A|_E=2l$, $F|_E=0$. 

\item $E\simeq\mathbb{P}^2$, $N_{E/X}=\mathcal{O}_{\mathbb{P}^2}(-2)$. In this case, we have $A\cdot l=1$ and the adjunction formula gives
\[
\mathcal{O}_E(A|_E) \otimes\mathcal{O}_E(2F|_E) = \mathcal{O}_E(1).
\]
As $A\cdot E$ is a non-zero effective cycle, we deduce that $F|_E=0$ and $A|_E=\mathcal{O}_E(1)$.
\end{enumerate}
We first show that $E$ cannot be contracted to a point. Suppose that we are in one of the above cases, then $F\cdot E=0$ and $A\cdot E$ is a non-zero effective cycle of rational curves. On the other hand, $A$ is a ruled surface over an elliptic curve, which implies that $E|_A$ consists of some fibers on $A$. But $F|_A$ is an elliptic curve which is a section, hence $E|_A\cdot F|_A >0$. This contradicts the fact that $F\cdot E=0$.

\item If $\varphi$ contracts $E$ to a smooth curve $C\subset Y$ of genus $g$, then $E\simeq\mathbb{P}(N_{C/Y}^*)$. Let $V=N_{C/Y}^*\otimes L$ with $L\in \Pic(C)$ be the normalization of the conormal bundle \cite[Chapter V, Proposition 2.8]{MR0463157}. Then  $N_{E/X}=-C_1+\mu l$ where $C_1$ is the minimal section satisfying $C_1^2=c_1(V)\eqqcolon -d$ and $\mu\coloneqq \dg L$. 

In this case, one has $-K_X\cdot l=1$, $F\cdot l=0$ and $A\cdot l=1$. Hence $F|_E=bl$ with $b\geq 0$ and the adjunction formula gives
\[
-K_E=(A+2F)|_E-E|_E,
\]
i.e. $A|_E=C_1 + (d+\mu+2(1-g-b))l$.

Since $F$ (resp. $A$) is a $\mathbb{P}^1$-bundle over a smooth elliptic curve, we deduce that the effective cycle $F\cdot E$ (resp. $A\cdot E$) does not contain the curve $l$ otherwise $l$ moves on the surface $F$ (resp. $A$). Therefore, $F\cdot E=0$ and $A\cdot E=C_1+(d+\mu+2(1-g))l$ is a section of $\varphi|_E\colon E\to C$. In particular, $\varphi(F)\simeq F$ as $E\cdot F=0$, and all the curves $l$ meet $A$ transversally in one point which implies
that $\varphi|_A$ is an isomorphism.

Now by the same argument as in the case (1), we deduce that the integral curve $A\cdot E$ cannot be a rational curve. Hence $C$ is of genus $g>0$. By \cite[Proposition 3.3]{MR1255695}, $-K_Y$ is again nef.

\end{enumerate}

Hence we have the following proposition:
\begin{prop}\label{biratfixedpart}
In Setup \ref{setupfixed}, let $\varphi\colon X\to Y$ be an extremal  contraction. Then $\varphi$ is the blow-up of a smooth curve $C$ of positive genus in the smooth threefold $Y$ with $-K_Y$ nef, $\kappa(-K_Y)=1$, $\nu(-K_Y)=2$. Let $E$ be the exceptional divisor of $\varphi$. Then one of the following two cases occurs:
\begin{enumerate}[label=\normalfont(\arabic*)]
\item $E=A$ and we have $|-K_Y|=|2G|$ with $G\coloneqq\varphi(F)\simeq F$. Furthermore, the blown up curve $C$ is a smooth elliptic curve contained in the base locus of $|G|$.
\item $E\neq A$ and $E\cdot F=0$. We have $|-K_Y|=A_Y+|2F_Y|$ where $A_Y\coloneqq\varphi(A)\simeq A$, $F_Y\coloneqq\varphi(F)\simeq F$ and $F_Y^2$ is a smooth elliptic curve. In particular, $Y$ satisfies again Setup \ref{setupfixed}.
\end{enumerate}
\end{prop}
\begin{proof}
\begin{enumerate}[label=(\arabic*)]
\item It remains to prove the last assumption of the first case.
Since $|-K_X|=A+|2H|$, one has $h^0(X,\mathcal{O}_X(H))=2$ by Corollary \ref{mH}. 

Now consider the threefold $Y$, since the anticanonical system $|-K_Y|=|2G|$ has no fixed part and again $-K_Y$ is nef with $n(-K_Y)=3$, $\nu(-K_Y)=2$, one has $h^0(Y,\mathcal{O}_Y(G))=2$ by the Corollary \ref{mH}.

Since $F$ is the strict transform of $G$ by $\varphi$, we deduce from $h^0(X,\mathcal{O}_X(F))=h^0(Y,\mathcal{O}_Y(G))$ that the blown up elliptic curve $C$ must be contained in the base locus of $|G|$.
\item 
Since $E\cdot F=0$, we have $\varphi^*(F_Y)=F$. We deduce that $F_Y$ is nef as $\varphi^*(F_Y)=F$ is nef.

We first show that $-K_Y$ is not semi-ample, which implies $\kappa(-K_Y)=1$ and $\nu(-K_Y)=2$.

Since $F^2$ is a non-zero effective $1$-cycle and $E\cdot F=0$, we deduce that $F_Y^2=\varphi(F)^2$ is also a non-zero effective $1$-cycle. Since $F_Y$ moves, $A_Y\cdot F_Y$ is an effective $1$-cycle. By the adjunction formula, we get
\[
-K_{F_Y}=(-K_Y-F_Y)|_{F_Y}=(A_Y+F_Y)|_{F_Y}
\]
and thus \[-K_Y|_{F_Y}=-K_{F_Y}+F_Y|_{F_Y}\]
is a non-zero effective divisor on $F_Y$ such that $-K_Y|_{F_Y}\leq -2K_{F_Y},$ i.e. $h^0(F_Y,\mathcal{O}_{F_Y}(-2K_{F_Y}-(-K_Y)|_{F_Y}))>0$.

Suppose by contradiction that $-K_Y$ is semi-ample, then $|-mK_Y|$ is base-point-free for $m\gg 0$. Hence its restriction $|-mK_Y|_{F_Y}|$ to $F_Y$ is also base-point-free.
On the other hand, since $F_Y\simeq F$, we have $\kappa(F_Y,-K_{F_Y})=0$. Hence
\[
1\leq h^0(F_Y,\mathcal{O}_{F_Y}(-mK_Y|_{F_Y}))\leq h^0(F_Y,\mathcal{O}_{F_Y}(-2mK_{F_Y}))=1.
\]
Therefore, the linear system $|-mK_Y|_{F_Y}|$ is fixed, which contradicts the fact that $|-mK_Y|_{F_Y}|$ is base-point-free.

Now we show that the anticanonical system $|-K_Y|$ has a fixed part. Since $F_Y$ is mobile, it is then clear that $A_Y$ is the fixed part of $|-K_Y|$.

Suppose by contradiction that $|-K_Y|$ has no fixed part, then $-K_Y$ has index two by the Theorem \ref{mainmobile}. As $-K_Y=A_Y+2F_Y$, this implies $A_Y=2L$ for some $L\in \Pic(Y)$. Hence
\[
-K_{F_Y}=(A_Y+F_Y)|_{F_Y}=(2L+F_Y)|_{F_Y}.
\]
Since $F_Y\simeq F$, $F_Y$ is a $\mathbb{P}^1$-bundle over a smooth elliptic curve such that $-K_{F_Y}\cdot f=2$ where $f$ is a fiber. Since $F^2$ is a smooth elliptic curve (a section of the $\mathbb{P}^1$-bundle $F$) and $E\cdot F=0$, we deduce that $F_Y^2$ is also a smooth elliptic curve (a section of the $\mathbb{P}^1$-bundle $F_Y$) and thus \[F_Y|_{F_Y}\cdot f=1.\]
This implies $2L|_{F_Y}\cdot f=1$, which contradicts the fact that $L|_{F_Y}$ is a Cartier divisor.
\end{enumerate}
 \end{proof}

\begin{rem}\label{elliptic}
In the setting of Proposition \ref{biratfixedpart} (2), we deduce by the same proposition that there exists a finite sequence \[X=X_0\overset{\varphi_1}\rightarrow X_1\overset{\varphi_2}\rightarrow\cdots\overset{\varphi_k}\rightarrow X_k\] where
\begin{itemize}
    \item $\varphi_i$ is a blow-up along a smooth curve $C_i$ of positive genus;
    \item $X_i$ satisfies again the setup \ref{setupfixed};
    \item $X_k$ has a birational extremal contraction which contracts the fixed part $A_k$ of $|-K_{X_k}|$.
\end{itemize}
Furthermore, the curve $C_i$ is contained in $A_i$, where $A_i$ is the fixed part of $|-K_{X_i}|$ which is a $\mathbb{P}^1$-bundle over a smooth elliptic curve $D_i$. Then $C_k$ is an elliptic curve and $k=1$.
\end{rem}
\begin{proof}
For $1\leq i\leq k$, let $g_i$ be a fiber of the $\mathbb{P}^1$-bundle $A_i$. 

Since $C_i$ has positive genus and it is contained in the $\mathbb{P}^1$-bundle $A_i$, it must be surjective to the curve $D_i$. Let $\alpha_i$ be the degree of $C_i$ onto the elliptic curve $D_i$. Then $g_i$ meets $C_i$ at $\alpha_i$ point(s). Hence in $X_{i-1}$, we have
\[E_{i-1}\cdot g_{i-1}=\alpha_i\] where $E_{i-1}$ is the exceptional divisor of $\varphi_i$ and $g_{i-1}$ is the strict transform of $g_i$. Therefore,
\[
-K_{X_{i-1}}\cdot g_{i-1}=\varphi_k^*(-K_{X_i})\cdot g_{i-1}-E_{i-1}\cdot g_{i-1}=-K_{X_i}\cdot g_i-\alpha_i.
\]
Since $-K_{X_{i-1}}$ is nef, we deduce that $-K_{X_i}\cdot g_i-\alpha_i\geq 0$ and thus $-K_{X_i}\cdot g_i\geq 1$.

For $i=k$, since $A_k$ is the exceptional divisor of an extremal contraction, we have $-K_{X_k}\cdot g_k=1$. Hence 
$\alpha_k=1$ (which implies $C_k\simeq D_k$ is a smooth elliptic curve) and $-K_{X_{k-1}}\cdot g_{k-1}=0$ (which implies $k=1$).
\end{proof}

\begin{proof}[Proof of Proposition \ref{mmp}.] It follows from Proposition \ref{biratfixedpart} and Remark \ref{elliptic}.
\end{proof}

\begin{proof}[Proof of Theorem \ref{mainfixed}.]
Let $|B|$ be the mobile part of the anticanonical system $|-K_X|$.
By Proposition \ref{nefmobile}, there exists a finite sequence of flops $\psi\colon X\dashrightarrow X'$ such that $-K_{X'}$ is nef and the mobile part $|B'|$ of $|-K_{X'}|$ is nef. 

Now we consider the case when $B$ is nef and suppose by contradiction that $B^2$ is a non-zero effective $1$-cycle. Then by Proposition \ref{mmp}, there exists a finite sequence \[X=X_0\overset{\varphi_1}\rightarrow X_1\overset{\varphi_2}\rightarrow\cdots \overset{\varphi_k}\rightarrow X_k\overset{\varphi_{k+1}}\rightarrow Y\] with $k=0$ or $1$, where $\varphi_i$ is a blow-up along a smooth elliptic curve and
$Y$ is one of the cases described in Theorem \ref{mainmobile} with $|-K_Y|=|2G|$. Moreover, a general member $D\in |G|$ is isomorphic to $F$, where $F$ is a general member in $|H|$. Hence $D$ is a $\mathbb{P}^1$-bundle over a smooth elliptic curve as described in Lemma \ref{classification}, $2.(B)$. 

On the other hand, $D$ is in one of the following cases:
\begin{enumerate}
    \item If $Y$ is a del Pezzo fibration: $\phi\colon Y\to\mathbb{P}^1$, then by Remark \ref{remdP}, $\phi\colon D\to\mathbb{P}^1$ induces a fibration on $D$ with general fiber isomorphic to $\mathbb{P}^1$ or two $\mathbb{P}^1$'s intersecting at one point.
    \item If $Y=\mathbb{P}_S(\mathcal{E})$ is a $\mathbb{P}^1$-bundle over a smooth rational surface $S$, where $\mathcal{E}$ is a nef rank-$2$ vector bundle on $S$ given by an extension
\[
0\to\mathcal{O}_S\to \mathcal{E}\to\mathcal{I}_Z\otimes \mathcal{O}_S(-K_S)\to 0
\]
with $\mathcal{I}_Z$ the ideal sheaf of $c_2(\mathcal{E})$ points on $S$,
then by Remark \ref{remconic}, $D=\Bl_Z(S)$ is a rational surface.
    \item If $Y$ has a birational extremal contraction, then $D$ contains a $(-1)$-curve by Remark \ref{rembirat}.
\end{enumerate}
Hence $D\not\simeq F$, which gives a contradiction.
\end{proof}

\newpage
\bibliographystyle{alpha}
  \bibliography{bibliography.bib}

\Affilfont{\small{Z\textsc{hixin} X\textsc{ie}, U\textsc{niversit\"{a}t} \textsc{des} S\textsc{aarlandes}, Saarbr\"{u}cken, Germany}}

\textit{Email address:}
\href{mailto:xie@math.uni-sb.de}{xie@math.uni-sb.de}
\end{document}